\documentclass[12pt, reqno]{amsart}
\usepackage{amsmath}
\usepackage{amssymb}
 2

\theoremstyle{remark}

\newcommand{\Z}{{\mathbb Z}}
\newcommand{\Q}{{\mathbb Q}}

\newcommand{\C}{{\mathbb C}}

\begin{document}

\title{The Strong Chowla-Milnor spaces and a conjecture of Gun, Murty and Rath}
\author{Tapas Chatterjee}
\address[Tapas Chatterjee]
        {The Institute of Mathematical Sciences,
         C.I.T Campus, Taramani,
         Chennai 600 113 
         India}
\email[Tapas Chatterjee]{tapasc@imsc.res.in}

\maketitle

\begin{abstract}
In a recent work, Gun, Murty and Rath formulated the Strong Chowla-Milnor conjecture and defined the 
Strong Chowla-Milnor space. In this paper, we prove a non-trivial lower bound for the dimension of 
these spaces. We also  obtain a conditional improvement of this lower bound and noted that an 
unconditional improvement of this lower bound will lead to irrationality of both $\zeta(k)$ and  
$\zeta(k)/ \pi^k$ for all odd positive integers $k>1$. Following Gun, Murty and Rath, we define
generalized Zagier spaces $V_p(K)$ for multiple zeta values over a number field $K$. We prove that the
dimension of $V_{4d+2}(K)$ for $d\geq 1$, is at least 2, assuming a conjecture of Gun, Murty and Rath. 
\end{abstract}

\section{{\bf Introduction}}

\bigskip

Throughout the paper, we consider $s\in \C$ with $\Re(s)>1$, unless otherwise stated. The Riemann zeta function is 
defined as
\begin{equation*}
 \zeta(s)=\sum_{\substack{n=1}}^\infty \frac{1}{n^s}, ~\text {for} ~\Re(s)>1.
\end{equation*}
The Riemann zeta function can be extended holomorphically to the whole complex plane except at $s=1$, where 
it has a simple pole with residue 1. Hurwitz studied the function, now called the Hurwitz zeta function,
which is defined as
\begin{equation*}
\zeta(s,x)=\sum_{n=0}^\infty \frac{1}{(n+x)^s}
\end{equation*}
where $0<x\leq 1$ and $s\in \C$ with $\Re(s)>1$. He proved that $\zeta(s,x)$ can be extended holomorphically to the 
entire complex plane except at $s=1$, where it has a simple pole with residue 1. Note that $\zeta(s,1)=\zeta(s)$, the 
classical Riemann zeta function.

For a periodic arithmetic function $f$ with period $q>1$ and $s\in\C$, the $L$-function associated to $f$ is defined as 
\begin{equation*}
  L(s,f)=\sum_{n=1}^\infty\frac{f(n)}{n^s}. 
\end{equation*}
 
 Since $f$ is periodic with period $q$, we have
\begin{equation*}
L(s,f)=q^{-s}\sum_{a=1}^q f(a)\zeta(s,a/q), 
\end{equation*}
Hence $L(s,f)$ extends holomorphically to the whole complex plane with a possible simple pole at 
$s=1$ with residue $q^{-1}\sum_{a=1}^{q} f(a)$.

For an algebraic valued periodic function $f$ with period $q$, the transcendental nature of $L(1,f)$, whenever it exists, is 
discussed in the paper \cite{MS1}. We got our motivation for this work from a conjecture formulated by Gun, Murty and Rath in
\cite{GMR} which is a generalization of the Chowla-Milnor conjecture stated below.

\bigskip
\noindent
\textbf{Conjecture} 1 (Chowla-Milnor).  Let $k>1 ~{\rm and} ~q>1$, be integers. Then the following $\varphi(q)$ real numbers 
are linearly independent over $\Q$:
\begin{equation*}
\zeta(k,a/q) ~ {\rm with} ~ 1\leq a < q, ~(a,q)=1.
\end{equation*}

The authors in \cite{GMR} generalized the above conjecture in the following manner:

\smallskip
\noindent
\textbf{Conjecture} 2 (Strong Chowla-Milnor). For any integers $k>1 ~{\rm and} ~q>1$, the following $\varphi(q)+1$ real numbers
\begin{equation*}
 1,\zeta(k,a/q) ~ {\rm with} ~ 1\leq a < q, ~(a,q)=1
\end{equation*}
are $\Q$-linearly independent. \\

For a rational valued periodic function $f$ with period $q$ and satisfying $f(a)=0$ for $1<(a,q)<q$, the Strong Chowla-Milnor conjecture 
is equivalent to the irrationality of $L(k,f)$, unless 
\begin{equation*}
f(a)=-\frac{f(q)q^{-k}}{\underset{p \in {\rm P}, \atop p|q}{\prod}(1-p^{-k})}. 
\end{equation*}
for $1\leq a<q, (a,q)=1$. Here $P$ denote the set of primes.\\

In order to see that the above statement is equivalent to the Strong Chowla-Milnor conjecture, we consider the identity
\begin{equation*}
\zeta(k)\underset{p \in {\rm P}, \atop  p|q}{\prod}(1-p^{-k})=
q^{-k}\sum_{\substack{a=1\\ (a,q)=1}}^{q-1}\zeta(k,a/q).
\end{equation*}
Substituting this in the above expression for $L(s,f)$ and using $f(a)=0$ for $1<(a,q)<q$, we get
\begin{equation}\label{1}
 L(k,f)=q^{-k}\sum_{\substack{a=1 \\ (a,q)=1}}^{q-1}\left[f(a)+\frac{f(q)q^{-k}}{\underset{p|q}\prod(1-p^{-k})}\right]\zeta(k,a/q). 
\end{equation}
If $L(k,f)$ is rational, then the above equation shows that $1, \zeta(k,a/q)$ for $1\leq a<q, (a,q)=1$ are linearly
dependent over the rationals since $f$ is a rational-valued function. Conversely, if $1, \zeta(k,a/q)$ with 
$1\leq a<q, (a,q)=1$ are linearly dependent over the rationals, then there are rational numbers $c_0, c_a, 1\leq a<q, (a,q)=1$,
not all zero, such that
\begin{equation*}
c_0+\sum_{\substack{a=1 \\ (a,q)=1}}^{q-1}c_a\zeta(k,a/q)=0. 
\end{equation*}
Now define the following rational-valued periodic function $f$ with period $q$. Set $f(a)=0$ for $1<a\leq q, (a,q)>1$ and 
$f(a)=c_a$ for $1\leq a<q, (a,q)=1$. Then, our identity shows that $q^kL(k,f)=-c_0$ so that $L(k,f)$ is rational. This 
shows that the above statement is equivalent to the Strong Chowla-Milnor conjecture.

\smallskip
\noindent
\textbf{Definition} 3. For an integer $m\geq 1$ and complex numbers $z$ with $|z|\leq 1$, the polylogarithm function
$Li_m(z)$ is defined by
\begin{equation*}
Li_m(z)=\sum_{n=1}^\infty \frac{z^n}{n^m}. 
\end{equation*}
Note that for $m=1$, the above series is equal to $-\log(1-z)$ for $|z|<1$. In a recent work \cite{GMR} Gun, Murty and
Rath formulated the Polylog conjecture. Here we have following strong conjecture about polylogarithms, which generalize 
Baker's theorem.

\smallskip
\noindent
\textbf{Conjecture} 4 (Strong Polylog Conjecture). Suppose $\alpha_1,...,\alpha_n$ are algebraic numbers with $|\alpha_i|\leq 1$
for $1 \leq i\leq n$, such that $Li_m(\alpha_1),..., Li_m(\alpha_n)$ are linearly independent over $\Q$.
Then $1, Li_m(\alpha_1),..., Li_m(\alpha_n)$ are linearly independent over $\overline{\Q}$.

Clearly for $m=1$ the above conjecture reduces to a consequence of Baker's theorem about linear forms in logarithms. In section 2, 
we establish a link between the Strong Polylog  conjecture and the Strong Chowla-Milnor conjecture in a form of 
a theorem stated as follows. 

\bigskip
\noindent
\textbf{Theorem} 5. The Strong Polylog conjecture implies the Strong Chowla-Milnor conjecture for all $q>1$ and $k>1$.

\bigskip
\noindent
\textbf{Definition} 6. For any integer $k>1$ and $q\geq 2$, define the Strong Chowla-Milnor space $\widehat{V}_k(q)$ by
\begin{equation*}
\widehat{V}_k(q) :=\Q- {\rm span ~of} ~\{1,\zeta(k,a/q):1\leq a < q, ~(a,q)=1\}. 
\end{equation*}

\smallskip
In a recent work \cite{GMR}, authors have shown that for an odd integer $k>1$ and two co-prime integers $q,r >2$, either 
\begin{equation*}
{\rm dim}_\Q V_k(q)\geq \frac{\varphi(q)}{2}+1
\end{equation*}
 or
\begin{equation*}
 {\rm dim}_\Q V_k(r)\geq \frac{\varphi(r)}{2}+1.
\end{equation*}
where the Chowla-Milnor space $V_k(q)$ is defined by 
\begin{equation*}
V_k(q) := \Q-{ \rm span ~ of~} \{\zeta(k,a/q):~ 1\leq a <q, ~(a,q)=1\}.
\end{equation*}

One can ask a similar type of question for the Strong Chowla-Milnor space i.e.,
for an odd integer $k>1$ and two co-prime integers $q,r >2$, whether or not either 
\begin{equation*}
{\rm dim}_{\Q}\widehat{V}_k(q)\geq \frac{\varphi(q)}{2}+2
\end{equation*}
 or
\begin{equation*}
 {\rm dim}_{\Q}\widehat{V}_k(r)\geq \frac{\varphi(r)}{2}+2.
\end{equation*}

But if the above statement is true, then clearly the Strong Chowla-Milnor conjecture is true for either
$q=3$ or $q=4$ i.e. either ${\rm dim}_{\Q}\widehat{V}_k(3)=3$ or ${\rm dim}_{\Q}\widehat{V}_k(4)=3$. Then 
from the proposition 7 of \cite{GMR} we get that $\zeta(k)$ is irrational for all odd $k>1$. In 
general we do not know for all odd integers $k$ whether $\zeta(k)$ is irrational or not. It is known,
thanks to Apery, that $\zeta(3)$ is irrational. On the other hand by a theorem of K. Ball and T. Rivoal 
(see \cite{TR} and \cite{BR}), it is known that $\zeta(k)$ is irrational for infinitely many odd $k>1$. 
In section 3, we will prove the following theorem.

\smallskip
\noindent
\textbf{Theorem} 7. Let $k>1$ be an odd integer with $\zeta(k)$ irrational and $q,r >2$ be two co-prime integers. 
Then either 
\begin{equation*}
{\rm dim}_{\Q}\widehat{V}_k(q)\geq \frac{\varphi(q)}{2}+2
\end{equation*}
 or
\begin{equation*}
 {\rm dim}_{\Q}\widehat{V}_k(r)\geq \frac{\varphi(r)}{2}+2.
\end{equation*}

\bigskip
In a recent work \cite{GMR1}, Gun, Murty and Rath formulated a variant of the Chowla-Milnor conjecture which is the following.

\smallskip
\noindent
\textbf{Conjecture} 8 (Gun, Murty and Rath). Let $K$ be a number field and  $k>1,q\geq 2$ be integers such that $K\cap \Q(\zeta_q)=\Q$. Then the 
$\varphi(q)$ real numbers $\zeta(k,a/q)$, $1\leq a < q , ~(a,q)=1$ are linearly independent over $K$.

\smallskip
In section 4, we investigate the values of $L(k,\chi)$, as $\chi$ ranges over non-trivial primitive 
Dirichlet characters mod $q$, over a certain family of algebraic number fields assuming the above conjecture.

\bigskip
\noindent
\textbf{Theorem} 9. Let $K$ be an algebraic number field and $K_1=K(e^{2\pi i/\varphi(q)})$. Suppose that $K_1\cap \Q(\zeta_q)=\Q$.
Then the values $L(k,\chi)$, as $\chi$ ranges over non-trivial primitive Dirichlet characters mod $q$, are linearly independent 
over $K_1$ for all $k\geq 1$, if the Gun, Murty and Rath conjecture is true.

\smallskip
In fact for $k>1$ the above theorem is true for any primitive Dirichlet  characters mod $q$, i.e. one can include the principal
character mod $q$.

\smallskip
In section 5, we establish a link between the Strong Chowla-Milnor conjecture and the multiple zeta values (MZVs).

\bigskip
\noindent
\textbf{Definition} 10. Let $k,s_1,...,s_k$ be positive integers with $s_1>1$. Then the \textit{multiple zeta values} (MZVs) are
defined as 
\begin{equation*}
\zeta(s_1,...,s_k)=\sum_{n_1>...>n_k\geq 1}\frac{1}{n_1^{s_1}...n_k^{s_k}}.
\end{equation*}

Clearly $k=1$ gives the classical Riemann zeta function. The sum $s_1+...+s_k$ is called the weight of the multiple zeta value
$\zeta(s_1,...,s_k)$ while $k$ is called the length of $\zeta(s_1,...,s_k)$.

\smallskip
\noindent
\textbf{Definition} 11. Let $K$ be a number field such that $K\cap \Q(\zeta_q)=\Q$. For any integer $p>1$, we define the 
generalized Zagier space as the $K$-linear space $V_p(K)$ defined by
\begin{equation*}
V_p(K)=K- {\rm span ~of} ~\{\zeta(s_1,...,s_k)| s_1+...+s_k=p\}. 
\end{equation*}

 In \cite{GMR}, authors have shown that the dimension of Zagier spaces $W_{4d+2}$ is at least 2 for all $d\geq 1$ assuming
the Chowla-Milnor conjecture, where the Zagier space is defined by
\begin{equation*}
W_p=\Q- {\rm span ~of} ~\{\zeta(s_1,...,s_k)| s_1+...+s_k=p\}. 
\end{equation*}
In this section we will prove the following theorem analogous to the theorem 3 in \cite{GMR}.

\smallskip
\noindent
\textbf{Theorem} 12. Let $d$ be a positive integer. Then the Gun, Murty and Rath conjecture implies 
\begin{equation*}
{\rm dim}_K V_{4d+2}(K)\geq 2. 
\end{equation*}

\section{\bf Proof of Theorem 5}

\smallskip
For the proof of Theorem 5, we shall need the following lemma (see S. Lang \cite{SL}, p.548).

\smallskip
\noindent
\textbf{Lemma} 13. Let $G$ be any finite abelian group of order $n$ and $F:G\rightarrow \C$ be any complex-valued function
on $G$. The determinant of the $n\times n$ matrix given by ($F(xy^{-1})$) as $x, y$ range over the group elements is 
called the {\it Dedekind determinant} and is equal to
\begin{equation*}
 \prod_\chi\left(\sum_{x\in G} \chi(x)F(x)\right),
\end{equation*}
where the product is over all characters $\chi$ of $G$. 

\smallskip
\noindent
\textbf{Proof of Theorem 5}. Let $k>1$ and $q>1$. Let $f$ be a rational valued period function with period $q$
satisfying $f(a)=0$ for $1<(a,q)<q$. Suppose that $L(k,f)=r$ is a rational number.
Then we have
\begin{equation*}
 L(k,f)=\sum_{n=1}^\infty \frac{f(n)}{n^k}=r
\end{equation*}

As $f$ is periodic function, we have the Fourier transformation of $f$ given by 
\begin{equation*}
\hat{f}(n)=\frac{1}{q}\sum_{a=1}^q f(a)\zeta^{-an}_q 
\end{equation*}
where $\zeta_q=e^{\frac{2\pi i}{q}}$ and hence we have the Fourier inversion formula
\begin{equation*}
f(n)=\sum_{a=1}^q \hat{f}(a)\zeta^{an}_q. 
\end{equation*}
Then we have
\begin{equation*}
 L(k,f)=\sum_{n=1}^\infty \frac{f(n)}{n^k}=\sum_{n=1}^\infty\frac{1}{n^k} \sum_{a=1}^q\hat{f}(a)\zeta^{an}_q=r
\end{equation*}
and hence 
\begin{equation*}
 \sum_{a=1}^q\hat{f}(a)Li_k(\zeta_q^a)-r=0.
\end{equation*}

Let $Li_k(\alpha_1),...,Li_k(\alpha_t)$ be a maximal linearly independent subset of
\begin{equation*}
 \{Li_k(\zeta_q^a)| 1\leq a < q\}  
\end{equation*}
over $\Q$.

Then \begin{equation*}
      Li_k(\zeta_q^a)=\sum_{b=1}^t C_{ab} Li_k(\alpha_b) 
     \end{equation*}
 
for some $C_{ab}\in \Q$. So we have 
\begin{equation*}
\sum_{b=1}^t \beta_b Li_k(\alpha_b)-r=0 
\end{equation*}
where
\begin{equation*}
 \beta_b=\sum_{a=1}^q \hat{f}(a)C_{ab}.
\end{equation*}

Since $f$ is rational valued, $\hat{f}$ is algebraic valued. So by the Strong Polylog conjecture, we have 
\begin{equation*}
r=0 ~{\rm and} ~ \beta_b=\sum_{a=1}^q \hat{f}(a)C_{ab}=0, ~1\leq b\leq t.
\end{equation*}

Now for any automorphism $\sigma$ of the field $\overline{\Q}$ over $\Q$, we have
\begin{equation*}
\sum_{a=1}^q \sigma(\hat{f}(a))C_{ab}=0, ~1\leq b\leq t,
\end{equation*}
and hence 
\begin{equation*}
\sum_{a=1}^q \sigma(\hat{f}(a))Li_k(\zeta_q^a)=0.
\end{equation*}
In particular, if for $1\leq h<q, (h,q)=1$, $\sigma_h$ is the element of the Galois group of $\Q(\zeta_q)$
over $\Q$ such that 
\begin{equation*}
\sigma_h(\zeta_q)=\zeta_q^h, 
\end{equation*}
then we have,
\begin{equation*}
 \sigma_h(\hat{f}(n))=\hat{f_h}(n)
\end{equation*}
where
\begin{equation*}
f_h(n)=f(nh^{-1}). 
\end{equation*}
Thus, we have
\begin{eqnarray*}
 L(k,f_h)&=&\sum_{n=1}^\infty \frac{f_h(n)}{n^k} \\ 
&=&\sum_{a=1}^q\hat{f_h}(a)Li_k(\zeta_q^a) \\
&=&\sum_{a=1}^q \sigma_h(\hat{f}(a))Li_k(\zeta_q^a)=0
\end{eqnarray*}
for all $1\leq h<q, (h,q)=1$. Thus by equation \eqref{1}, we get
\begin{equation*}
L(k,f_h)=q^{-k}\sum_{\substack{a=1 \\ (a,q)=1}}^{q-1}\left[f_h(a)+\frac{f_h(q)q^{-k}}{\underset{p|q}\prod(1-p^{-k})}\right]\zeta(k,a/q)=0 
\end{equation*}
for all $1\leq h<q, (h,q)=1$.

Now, putting $ah^{-1}=b$ and noting that $f_h(q)=f(q)$, we have
\begin{equation}\label{2}
L(k,f_h)=q^{-k}\sum_{\substack{b=1 \\ (b,q)=1}}^{q-1}\left[f(b)+\frac{f(q)q^{-k}}{\underset{p|q}\prod(1-p^{-k})}\right]\zeta(k,bh/q)=0
\end{equation}
for all $1\leq h<q, (h,q)=1$. 

Thus we get a matrix equation with $M$ being the $\varphi(q)\times\varphi(q)$ matrix whose $(b,h)$-th entry is given
by $\zeta(k,bh/q)$. Then by the evaluation of the Dedekind determinant as in Lemma 13, we get
\begin{equation*}
Det(M)=\pm\prod_\chi q^k L(k,\chi)\neq 0. 
\end{equation*}
Thus the matrix $M$ is invertible and hence by the equation \eqref{2}, we have 
\begin{equation*}
f(a)+\frac{f(q)q^{-k}}{\underset{p|q}\prod(1-p^{-k})}=0, ~1\leq a< q, ~(a,q)=1 
\end{equation*}
and hence
\begin{equation*}
f(a)=-\frac{f(q)q^{-k}}{\underset{p|q}\prod(1-p^{-k})} 
\end{equation*}
for all $1\leq a< q, (a,q)=1$. This completes the proof of theorem 5. 

\section{\bf Dimension of Strong Chowla-Milnor Spaces}  

The following lemma 14, due to Okada~\cite{OK} about the linear independence of co-tangent values at rational 
arguments, plays a significant role in proving the theorem 7.

\smallskip
\noindent
{\bf Lemma} 14. Let $k$ and $q$ be positive integers with $k\ge 1$ and $q>2$. Let T be a set of $\varphi(q)/2$
representations mod $q$ such that the union $T\cup(-T)$ constitutes a complete set of co-prime residue classes 
mod $q$. Then the set of real numbers
\begin{equation*}
 \frac{d^{k-1}}{dz^{k-1}}\cot(\pi z)|_{z=a/q}, ~ a\in T
\end{equation*}
\textit{is linearly independent over $\Q$}.

We first start with a proposition.

\bigskip
\noindent
\textbf{Proposition} 15.  $ \frac{\varphi(q)}{2}+1\leq \dim_{\Q}\widehat{V}_k(q) \leq \varphi(q)+1 $.

\noindent
\textbf{Proof}. Clearly from the definition of the Strong Chowla-Milnor space $\dim_{\Q}\widehat{V}_k(q) \leq \varphi(q)+1$. 
Note that the space $\widehat{V}_k(q)$ is also spanned by the following sets of real numbers:
\begin{equation*}
1,\{\zeta(k,a/q)+\zeta(k,1-a/q)|(a,q)=1, ~ 1\leq a<q/2\},
\end{equation*}
\begin{equation*}
 \{\zeta(k,a/q)-\zeta(k,1-a/q)|(a,q)=1, ~ 1\leq a<q/2\}. 
\end{equation*}

Then we have the following (see \cite{MS} and \cite{GMR}, for instance) 
\begin{equation*}
\zeta(k,a/q)+(-1)^k \zeta(k,1-a/q)=\frac{(-1)^{k-1}}{(k-1)!}\frac{d^{k-1}}{dz^{k-1}} (\pi\cot\pi z)|_{z=a/q} {\rm ~for} ~k\geq 2
\end{equation*}
and by Okada(\cite{OK}), we get $$\zeta(k,a/q)+(-1)^k \zeta(k,1-a/q)$$ for $1\leq a< q/2, (a,q)=1$,  are linearly independent over $\Q$.
Again by induction, we have 
\begin{equation*} 
\frac{d^{k-1}}{dz^{k-1}}(\pi\cot\pi z)=\pi^k \times \Z {\rm ~linear ~ combination ~of} ~(\csc\pi z)^{2l}(\cot\pi z)^{k-2l} ,
\end{equation*}
for some non-negative integer $l$ and for an integer $k\geq 1$. Since $\csc\pi z$ and $\cot\pi z$ are algebraic at
rationals, we have all the numbers $\zeta(k,a/q)+(-1)^k \zeta(k,1-a/q)$ are transcendental for any $k$. Hence
 $\dim_{\Q}\widehat{V}_k(q) \geq \frac{\varphi(q)}{2}+1$.

\bigskip
\noindent
\textbf{Proof of Theorem 7}. 
Suppose not, then we have 
\begin{equation*}
 {\rm dim}_{\Q}\widehat{V}_k(q)= \frac{\varphi(q)}{2}+1.
\end{equation*}
This gives that the numbers 
\begin{equation*}
1, \zeta(k,a/q)-\zeta(k,1-a/q), ~ {\rm where} ~(a,q)=1, ~ 1\leq a<q/2
\end{equation*}
generate $\widehat{V}_k(q)$. 

Since $k$ is odd, then by Hecke \cite{HE} (See also paper 41 of E. Hecke, Mathematische Werke, Dritte Auflage, Vandenhoeck
und Rupertecht, Gottingen, 1983), we have 
\begin{equation}\label{3}
 \frac{\zeta(k,a/q)-\zeta(k,1-a/q)}{(2\pi i)^k}\in \Q(\zeta_q).
\end{equation}
Again we know that
\begin{equation*}
\zeta(k)\underset{p \in {\rm P}, \atop p|q}{\prod}(1-p^{-k}) ~=~
q^{-k}\sum_{\substack{a=1,\\(a,q)=1}}^{q-1}\zeta(k,a/q) \in \widehat{V}_k(q),
\end{equation*}
where $P$ be the set of primes. \\
Thus $\zeta(k)\in \widehat{V}_k(q)$ and hence we have
\begin{equation*}
\zeta(k)=q_1+\sum_{\substack{(a,q)=1 \\ 1\leq a<q/2}}\lambda_a[\zeta(k,a/q)-\zeta(k,1-a/q)] ~{\rm for ~some}~ q_1,~\lambda_a\in \Q 
\end{equation*}
so that
\begin{equation*}
\frac{\zeta(k)-q_1}{(2\pi i)^k}=\sum_{\substack{(a,q)=1\\1\leq a<q/2}}\frac{\lambda_a[\zeta(k,a/q)-\zeta(k,1-a/q)]}{(2\pi i)^k}. 
\end{equation*}

\noindent
Thus by \eqref{3}
\begin{equation*}
\frac{\zeta(k)-q_1}{i\pi^k}=a_1({\rm say})\in \Q(\zeta_q).
\end{equation*}
Similarly, if 
\begin{equation*}
 {\rm dim}_{\Q}\widehat{V}_k(r)= \frac{\varphi(r)}{2}+1,
\end{equation*}
then
\begin{equation*}
 \frac{\zeta(k)-q_2}{i\pi^k}=a_2 ({\rm say})\in \Q(\zeta_r),  ~~{\rm with}  ~q_2\in \Q.
\end{equation*}

So we have
\begin{equation*}
 a_1i\pi^k+q_1=a_2i\pi^k+q_2
\end{equation*}
which implies
\begin{equation*}
(a_1-a_2)i\pi^k =q_2-q_1 .
\end{equation*}
The L.H.S of the above equation is algebraic number times transcendental number hence transcendental and the 
R.H.S is a rational number. Hence we get that $q_1=q_2$ and $a_1=a_2$.

Thus we have
\begin{equation*}
 \frac{\zeta(k)-q_1}{i\pi^k}\in \Q(\zeta_q)\cap\Q(\zeta_r)=\Q.
\end{equation*}

Let 
\begin{equation*}
 \frac{\zeta(k)-q_1}{i\pi^k}=a \in\Q.
\end{equation*}
Since L.H.S of the above equation is purely imaginary and R.H.S is rational, we have
$a=0$ and $\zeta(k)=q_1$, a rational number. This is a contradiction to the irrationality of 
$\zeta(k)$. Thus either 
\begin{equation*}
{\rm dim}_{\Q}\widehat{V}_k(q)\geq \frac{\varphi(q)}{2}+2
\end{equation*}
 or
\begin{equation*}
 {\rm dim}_{\Q}\widehat{V}_k(r)\geq \frac{\varphi(r)}{2}+2.
\end{equation*}
This completes the proof of the theorem.

\bigskip
\noindent
\textbf{Proposition} 16. $2\leq \dim_{\overline{\Q}}\widehat{V}_k(q) \leq \frac{\varphi(q)}{2}+2$.

\noindent
\textbf{Proof}. We know that $\zeta(k,a/q)+(-1)^k \zeta(k,1-a/q)\in \pi^k\overline{\Q}^*$ for $1\leq a < q/2, (a,q)=1$.
As 1 and $\pi$ are linearly independent over $\overline{\Q}$, we have $\dim_{\overline{\Q}}\widehat{V}_k(q) \geq 2$.

As $\zeta(k,a/q)+(-1)^k \zeta(k,1-a/q)$ for $1\leq a < q/2, (a,q)=1$, are linearly dependent over $\overline{\Q}$, they contribute 
at most 1 in the dimension. Hence we have $\dim_{\overline{\Q}}\widehat{V}_k(q) \leq \frac{\varphi(q)}{2}+2$.

\section{\bf Proof of Theorem 9}

\smallskip
\noindent
\textbf{Proof}. The case $k=1$ is a theorem of M. Ram Murty and N. Saradha \cite{MS}. Let $k>1$ and assume that
\begin{equation*}
\sum_{\chi} c_\chi L(k,\chi)=0,  ~ c_\chi\in K_1 
\end{equation*}
where the summation is over all primitive Dirichlet characters mod $q$.

Again we know that 
\begin{equation*}
L(k,\chi)=q^{-k}\sum_{\substack{1\leq a< q \\ (a,q)=1}} \chi(a)\zeta(k,a/q). 
\end{equation*}

So from the above equation we get 
\begin{equation*}
q^{-k}\sum_{\chi} c_\chi \sum_{\substack{1\leq a< q \\ (a,q)=1}} \chi(a)\zeta(k,a/q)=0, 
\end{equation*}
 and hence we have
\begin{equation*}
 \sum_{\substack{1\leq a< q \\ (a,q)=1}} \zeta(k,a/q)\sum_{\chi} c_\chi \chi(a)=0. 
\end{equation*}
The values of $\chi$ lie in the field $K_1$ so that the sum $\underset{\chi}{\sum} c_\chi \chi(a)\in K_1$ which is disjoint from 
$\Q(\zeta_q)$. Hence using Gun, Murty and Rath conjecture, we get
\begin{equation*}
\sum_{\chi} c_\chi \chi(a)=0 
\end{equation*}
for all $a\in (\Z/qZ)^*$.
Then by the orthogonality of characters, we have $c_\chi=0$ for all $\chi$. This completes the proof. 

\section{\bf Dimension of generalized Zagier spaces}

Before we proceed to the proof of the above theorem 12, we review some basic facts about Kronecker symbols and associated Gauss sums.
If $\Delta$ be a fundamental discriminant then we have $\Delta\equiv 0, 1$ (mod 4). Then the Kronecker symbol $(\frac{\Delta}{n})$ 
is defined by the following relations:

(i) $(\frac{\Delta}{p})=0$ when $p|\Delta$ and $p$ prime, \\

(ii)$(\frac{\Delta}{2})=\left\{
                        \begin{array}{ll}
                         1 ~ when ~\Delta\equiv 1 ~(mod ~8), \\
                         -1 ~when ~\Delta\equiv 5 ~(mod ~8),
                        \end{array}\right. $

(iii)$(\frac{\Delta}{p})=(\frac{\Delta}{p})_L$, the Legendre symbol, when $p>2$,

(iv)$(\frac{\Delta}{-1})=\left\{
                        \begin{array}{ll}
                         1 ~ when ~\Delta >0, \\
                         -1 ~when ~\Delta <0,
                        \end{array}\right. $
  
(v)$(\frac{\Delta}{n})$ is completely multiplicative function of $n$.

In the proof of theorem 12, we use the following theorem( see \cite{MV}, p.297).

\smallskip
\noindent
\textbf{Theorem} 17. Let $\Delta$ be a fundamental discriminant. Then $\chi_\Delta (n)=(\frac{\Delta}{n})$ is a primitive quadratic character 
modulo $|\Delta|$.

Clearly from the definition $\chi_\Delta(n)$ is an odd character if $\Delta < 0$.

We first start with a lemma.

\smallskip
\noindent
\textbf{Lemma} 18. Let $K$ be an algebraic number field. Then $\left[\frac{\zeta(2d+1)}{\pi^{2d+1}}\right]^2 \notin K $ implies 
 ${\rm dim}_K V_{4d+2}(K)\geq 2$.

\smallskip
\noindent
\textbf{Proof}. Using the definition of multiple zeta values we get 
\begin{equation*}
\zeta(s_1)\zeta(s_2)=\zeta(s_1,s_2)+\zeta(s_2,s_1)+\zeta(s_1+s_2).
\end{equation*}
Thus we have
\begin{equation*}
\zeta(2d+1)^2=2\zeta(2d+1,2d+1)+\zeta(4d+2)
\end{equation*}
and
\begin{equation*}
\left[\frac{\zeta(2d+1)}{\pi^{2d+1}}\right]^2=2\frac{\zeta(2d+1,2d+1)}{\pi^{4d+2}}+\frac{\zeta(4d+2)}{\pi^{4d+2}}. 
\end{equation*}
Since 
\begin{equation*}
\left[\frac{\zeta(2d+1)}{\pi^{2d+1}}\right]^2 \notin K ~and ~\frac{\zeta(4d+2)}{\pi^{4d+2}}\in K 
\end{equation*}
it follows that $\zeta(2d+1,2d+1)$ is not in the $K$- span of $\zeta(4d+2)$ and hence the $K$-dimension of the space
$V_{4d+2}(K)\geq 2$.

\smallskip
\noindent
\textbf{Lemma} 19. Suppose the Gun, Murty and Rath conjecture is true. Then 
\begin{equation*}
 \left[\frac{\zeta(2d+1)}{\pi^{2d+1}}\right]^2 \notin K, 
\end{equation*}
for all $d\geq 1$.

\smallskip
\noindent
\textbf{Proof}. Let $\Delta<0$ be a fundamental discriminant. Then the Kronecker symbol $\chi_\Delta (n)=(\frac{\Delta}{n})$
is an odd, primitive, quadratic character modulo $|\Delta|$. Let $q=|\Delta|$. 

Now from the theory of Gauss sums we know that the Gauss sums $\tau(\chi_\Delta)$ associated with the Kronecker symbol $\chi_\Delta$
(see \cite{MV}, p.300) is given by
\begin{equation*}
 \tau(\chi_\Delta)=\sum_{a=1}^q \chi_\Delta(a)\zeta_q^a=i\surd q.
\end{equation*}
Again using the primitivity of $\chi_\Delta$, we have
\begin{equation*}
 \tau(\chi_\Delta, b)=\sum_{a=1}^q \chi_\Delta(a)\zeta_q^{ab}=\overline{\chi}_\Delta(b)i\surd q.
\end{equation*}
Since $\chi_\Delta$ is an odd character, we have
\begin{equation*}
\sum_{a=1}^{q/2}\chi_\Delta(a)(\zeta_q^{ab}-\zeta_q^{-ab})=\overline{\chi}_\Delta(b)i\surd q. 
\end{equation*}
 Let $B_l(x)$ be the $l$th Bernoulli polynomial. Multiplying both sides of the above equation by $B_{2d+1}(b/q)$
and taking sum over $b=1$ to $q$ we get,
\begin{equation*}
\sum_{a=1}^{q/2}\chi_\Delta(a)\sum_{b=1}^q(\zeta_q^{ab}-\zeta_q^{-ab})B_{2d+1}(b/q)=i\surd q\sum_{b=1}^q\overline{\chi}_\Delta(b)B_{2d+1}(b/q). 
\end{equation*}

Let $k=2d+1$. Then from proposition 1 of \cite{GMR}, we have 
\begin{equation*}
 \frac{\zeta(k,a/q)-\zeta(k,1-a/q)}{(2\pi i)^k}=\frac{q^{k-1}}{2k!}\sum_{b=1}^q(\zeta_q^{ab}-\zeta_q^{-ab})B_k(b/q)
\end{equation*}
for any $(a,q)=1$ and $1\leq a < q/2$. As $\chi_\Delta$ is a quadratic character we get that the number $i\surd q$ lies in the $K$-linear
space generated by the real numbers
\begin{equation*}
 \frac{\zeta(k,a/q)-\zeta(k,1-a/q)}{(2\pi i)^k}
\end{equation*}
with $(a,q)=1$ and $1\leq a < q/2$.
    
Again we know that 
\begin{equation*}
\zeta(k)\underset{p \in {\rm P}, \atop p|q}{\prod}(1-p^{-k})=
q^{-k}\sum_{\substack{a=1\\ (a,q)=1}}^{q-1}\zeta(k,a/q) 
\end{equation*}
where $P$ be the set of primes. So that 
\begin{equation*}
\zeta(k)\underset{p \in {\rm P}, \atop p|q}{\prod}(1-p^{-k})=
q^{-k}\sum_{\substack{a=1\\ (a,q)=1}}^{q/2}[\zeta(k,a/q)+\zeta(k,1-a/q)].  
\end{equation*}

Hence $\zeta(k)/(2\pi i)^k$ lies in the $K$-linear space generated by the real numbers
\begin{equation*}
\frac{\zeta(k,a/q)+\zeta(k,1-a/q)}{(2\pi i)^k} 
\end{equation*}
with $(a,q)=1$ and $1\leq q < q/2$.

Thus the Gun, Murty and Rath conjecture for the modulus $q$ implies that $i\surd q$ and $\zeta(k)/(2\pi i)^k$ lie in two disjoint
$K$-spaces. Hence for any such $q$, we have
\begin{equation*}
 \frac{\zeta(2d+1)}{\pi^{2d+1}\surd q}\notin K.
\end{equation*}
Thus, if the Gun, Murty and Rath conjecture true for all modulus, then 
\begin{equation*}
\left[\frac{\zeta(2d+1)}{\pi^{2d+1}}\right]^2 \notin K 
\end{equation*}
for all $d\geq 1$.

\smallskip
\noindent
\textbf{Proof of theorem} 12. Suppose the Gun, Murty and Rath conjecture is true. Then lemma 19 implies
\begin{equation*}
\left[\frac{\zeta(2d+1)}{\pi^{2d+1}}\right]^2 \notin K. 
\end{equation*}
 Hence from lemma 18, we get
\begin{equation*}
{\rm dim}_K V_{4d+2}(K)\geq 2. 
\end{equation*}
This completes the proof of theorem 12. 

\bigskip
\noindent
{\bf Acknowledgement.}
I would like to thank Sanoli Gun for many valuable suggestions. I would also like to thank the referee for several helpful comments.

\end{document}